\documentclass[11pt]{article}
	\usepackage[utf8]{inputenc}

\usepackage{amsmath,amssymb,amstext,dcolumn,epsfig,algorithmic,url}
\usepackage{geometry}
\usepackage{xcolor}

\usepackage{makecell}
\usepackage{arydshln}

\newtheorem{algorithm2}{Algorithm}
\newenvironment{algorithm}[1]{\begin{algorithm2}\rm{ #1}\begin{algorithmic}[1]}{\end{algorithmic}\end{algorithm2}}

\newcommand{\C}{{\mathbb C}}

\def\stril{\mathop{\mathsf{stril}}\nolimits}

\newtheorem{theorem}{Theorem}
\newtheorem{example}[theorem]{Example}
\newtheorem{lemma}[theorem]{Lemma}
\newtheorem{remark}[theorem]{Remark}

\newcommand {\proof} {\par{\it Proof}. \ignorespaces}
\newcommand {\eproof}
      {\space
        {\ \vbox{\hrule\hbox{\vrule height1.3ex\hskip0.8ex\vrule}\hrule}}
        \par}

\newcommand{\hide}[1]{}

\newcommand{\code}[1]{\texttt{#1}}

\catcode`@=11
\font\tenex=cmex10 
\newdimen\p@renwd
\setbox0=\hbox{\tenex B} \p@renwd=\wd0 
\def\bmat#1{\begingroup \m@th
  \setbox\z@\vbox{\def\cr{\crcr\noalign{\kern2\p@\global\let\cr\endline}}%
    \ialign{$##$\hfil\kern2\p@\kern\p@renwd&\thinspace\hfil$##$\hfil
      &&\quad\hfil$##$\hfil\crcr
      \omit\strut\hfil\crcr\noalign{\kern-\baselineskip}%
      #1\crcr\omit\strut\cr}}%
  \setbox\tw@\vbox{\unvcopy\z@\global\setbox\@ne\lastbox}%
  \setbox\tw@\hbox{\unhbox\@ne\unskip\global\setbox\@ne\lastbox}%
  \setbox\tw@\hbox{$\kern\wd\@ne\kern-\p@renwd\left[\kern-\wd\@ne
    \global\setbox\@ne\vbox{\box\@ne\kern2\p@}%
    \vcenter{\kern-\ht\@ne\unvbox\z@\kern-\baselineskip}\,\right]$}%
  \null\;\vbox{\kern\ht\@ne\box\tw@}\endgroup}
\catcode`@=12

\title{Iterative Refinement of Schur decompositions}
\author{Zvonimir Bujanovi\'c\footnote{University of Zagreb, Faculty of Science, Department of Mathematics, Croatia. {\tt zbujanov@math.hr}. Supported by the Croatian Science Foundation under grant IP-2019-04-6268.} \and Daniel Kressner\footnote{Institute of Mathematics, EPFL, Switzerland. {\tt daniel.kressner@epfl.ch}} \and Christian Schr\"oder\footnote{Freelancing numerical analyst, last academic position: Institut f\"ur Mathematik, TU Berlin, Germany. {\tt chris.schroeder@gmx.net}}}

\begin{document}

\maketitle

\begin{abstract}
The Schur decomposition of a square matrix $A$ is an important intermediate step of  state-of-the-art numerical algorithms for addressing eigenvalue problems, matrix functions, and matrix equations. This work is concerned with the following task: Compute a (more) \emph{accurate} Schur decomposition of $A$ from a given \emph{approximate} Schur decomposition. This task arises, for example, in the context of parameter-dependent eigenvalue problems and mixed precision computations. We have developed a Newton-like algorithm that requires the solution of a triangular matrix equation and an approximate orthogonalization step in every iteration. We prove local quadratic convergence for matrices with mutually distinct eigenvalues and observe fast convergence in practice. In a mixed low-high precision environment, our algorithm essentially reduces to only four high-precision matrix-matrix multiplications per iteration. When refining double to quadruple precision, it often needs only $3$--$4$ iterations, which reduces the time of computing a quadruple precision Schur decomposition by up to a factor of $10$--$20$.
\end{abstract}

\pagestyle{myheadings}
\thispagestyle{plain}
\markboth{D. KRESSNER, Z. BUJANOVI\'C, C. SCHR\"ODER}{Refinement of Schur decompositions}

\section{Introduction}

Given a matrix $A \in \C^{n\times n}$, a factorization of the form
\begin{equation} \label{eq:complexschur}
 T = Q^H A Q,
\end{equation}
with $Q \in \C^{n\times n}$ unitary and $T\in \C^{n\times n}$ upper triangular is called \emph{Schur decomposition} of $A$. 
This decomposition plays a central role in algorithms for solving eigenvalue problems, computing matrix functions, and solving matrix equations.
Note that the diagonal of $T$ contains the eigenvalues of $A$ and they can be reordered to appear there in any desirable order~\cite{Golub2013}.

In this work, we consider the following question. Suppose that we have an approximate Schur decomposition at our disposal, that is, $T$ is \emph{nearly} upper triangular and $Q$ is unitary or only \emph{nearly} unitary. How do we refine $Q,T$ to yield a more accurate Schur decomposition using an algorithm that is computationally and conceptually less demanding than computing the Schur decomposition of $A$ from scratch by applying, e.g., the QR algorithm~\cite{BraBM02a,BraBM02b}?

Partly motivated by the increased role of GPUs and TPUs in high-performance computing, there has been a revival of interest in exploiting the benefits of a mixed-precision environment in numerical computations; see~\cite{Abdelfattah2021} for a recent survey. Now, suppose that a Schur decomposition of $A$ has been computed, inexpensively, in a certain (low) machine precision but the target application demands for higher accuracy. For the computed factor $\hat Q$, the entries of $\hat Q^H A \hat Q$ below the diagonal and the entries of $\hat Q^H \hat Q - I$ are roughly on the level of unit roundoff in low precision. The refinement procedure discussed in this work produces a Schur decomposition that is accurate in high machine precision essentially at the cost of a few matrix multiplications in high machine precision. Examples for this setting of current interest include mixing half or single with double precision on specialized processors as well as mixing double with quadruple or even higher precision on standard CPUs. In both cases, the operations performed in high precision are significantly more costly and should therefore be limited to
the minimum. Other scenarios that could benefit from the refinement of Schur decompositions include the solution of parameter-dependent eigenvalue problems and continuation methods~\cite{BeyKT01,BinDF08}.

The goal of this paper is to develop an efficient Newton-like method for refining a Schur decomposition.
Closely related to the theme of this paper, the refinement of an individual eigenvector has been investigated
quite thoroughly, going back to the works of Wielandt~\cite{Wie44} and Wilkinson~\cite{Wil65}, and is usually addressed by some form of inverse iteration, see also~\cite{Ips97}. When several eigenvalues and/or eigenvectors are of interest, it is sensible to refine the invariant subspace (belonging to the eigenvalues of interest) as a whole rather than individual eigenvectors. Various variants of the Newton algorithm have been proposed for this purpose, see, e.g.,~\cite{AbsMSV02,BeyKT01,BinDF08,Cha84,Dem87d,LunE02}.
Often, these algorithms require the solution of a Sylvester equation at each iteration. 
While the developments presented in this work bear similarities, we are not aware that any of these existing methods would allow for refining a Schur decomposition or, equivalently, for simultaneously refining the entire flag of invariant subspaces associated with the Schur decomposition. In principle, Jacobi algorithms could be used for this purpose; see~\cite{Ebe62,Gre55,HueV03,Ste85} for examples. While such algorithms converge locally and asymptotically quadratically under certain conditions~\cite{Meh08}, their efficient implementation requires significant attention to low-level implementation aspects while the critical parts of our algorithm are entirely based on matrix multiplications. Improving an earlier method by Davies and Modi~\cite{Davies1986}, Ogita and Aishima~\cite{Ogita18,Ogita2019} recently proposed a Newton-type methods for refining spectral decompositions, addressing the task considered in this work when $A$ is real symmetric/complex Hermitian matrix. We will discuss similarities and differences with our algorithm in Section~\ref{sec:symmetrica}.
As a curiosity, we also point out a discussion by Kahan~\cite{Kahan2007}
on refining eigendecompositions for nonsymmetric matrices; see also Jahn's 1948 work~\cite{Jahn1948}.


\section{Algorithms}

To motivate the approach pursued in this work, let us consider a unitary matrix $\hat Q$ that nearly effects a Schur decomposition:
\begin{equation} \label{eq:lowup}
 \hat Q^H A \hat Q = T + E, \qquad \varepsilon := \|E\|_F \approx 0,
\end{equation}
where $T$ is upper triangular and $E$ is strictly lower triangular. We will write $\stril(\cdot)$ to denote the strictly lower triangular part of a matrix. A unitary matrix $Q_L$ that transforms $T+E$ to Schur form satisfies the two matrix equations
\begin{equation} \label{eq:nonlinearsystem}
 \stril(Q_L^H (T+E) Q_L) = 0, \quad Q_L^H Q_L = I.
\end{equation}
In the following we linearize these equations, analogous to existing first-order perturbation analyses of Schur decompositions~\cite{KonPC94,Sun92}.
Writing $Q_L = I + W$, the condition $Q_L^H Q_L = I$ becomes equivalent to
\begin{equation} \label{eq:orthW}
 W + W^H + W^H W = 0
\end{equation}
Ignoring the second order term $W^H W$, this means that $W$ is skew-Hermitian and can be written as
\begin{equation} \label{eq:W}
 W = L + D - L^H, \quad \text{with }\ L = \stril(W),
\end{equation}
where $D$ is diagonal with purely imaginary diagonal elements. Now, the first equation in~\eqref{eq:nonlinearsystem} becomes
\[
 \stril\left( (I+L+D-L^H)^H (T+E) (I+L+D-L^H) \right) = 0.
\]
Assuming $\|W\|_F = O(\varepsilon)$, where $\|\cdot\|_F$ denotes the Frobenius norm, and dropping all second order terms in $\varepsilon$, we arrive at the triangular matrix equation
\begin{equation} \label{eq:linearsystem}
 \stril(E - L T + TL) = 0.
\end{equation}
If the eigenvalues of $T$ are pairwise distinct, there is a unique strictly lower triangular matrix $L$ satisfying~\eqref{eq:linearsystem},
see~\cite{KonPC94,Sun92} and also Theorem~\ref{thm:solvability} below.
Note that the diagonal factor $D$ has disappeared in~\eqref{eq:linearsystem};
it can be chosen to contain arbitrary diagonal entries on the imaginary axis.
In the following, we choose $D = 0$.

In order to attain a unitary factor $Q_L$, the matrix $I + W = I+L-L^H$ needs to be orthogonalized; we will discuss different orthogonalization strategies in Section~\ref{sec:orthog}.

Algorithm~\ref{alg:template} summarizes the procedure outlined above. Note that we allow the input matrix $\hat Q$ to be non-unitary, which is taken care of by an optional orthogonalization step in line~\ref{line:optqr}.
\begin{algorithm}{Template for refining a Schur decomposition}\label{alg:template}
\REQUIRE Matrix $A \in \C^{n\times n}$ and (nearly) unitary matrix $\hat Q \in \C^{n\times n}$  such that $\hat Q^H A \hat Q$ is nearly upper triangular.
\ENSURE Unitary matrix $Q$ such that $Q^H A Q$ is upper triangular.

 \STATE Compute (approximately) unitary matrix $Q$ from orthogonalizing $\hat Q$.  \label{line:optqr}
 \REPEAT
   \STATE \label{line:extractT} Compute $\hat T \gets Q^H A Q$.
   \STATE Set $E\leftarrow \stril(\hat T)$, $T\leftarrow \hat T - E$.
   \STATE \label{line:solvebasic}Solve equation~(\ref{eq:linearsystem}) for strictly lower triangular $L$.
   \STATE \label{line:unitary} Replace $Q$ by (approximately) unitary matrix obtained from orthogonalizing $Q(I+L-L^H)$.
  \UNTIL{convergence}
\end{algorithm}

In the following sections, we will provide details on Algorithm~\ref{alg:template} and perform a convergence analysis.

\subsection{Solution of triangular matrix equation}
\label{sec:solveSylv}

Let us first consider the triangular matrix equation from Step~\ref{line:solvebasic} of Algorithm~\ref{alg:template}:
\begin{equation} \label{eq:linearsystemagain}
 \stril(TL - LT) = -E,
\end{equation}
where $T$ is upper triangular and $E$ as well as the desired solution $L$ are strictly lower triangular. For $1\le j<i\le n$, equating the $(i,j)$-entry of~\eqref{eq:linearsystemagain} yields
\begin{equation}\label{eq:aij}
 \sum_{k=i}^n t_{ik}\ell_{ki}-\sum_{k=1}^j \ell_{ik} t_{kj}=-e_{ij},
\end{equation}
which, assuming $t_{ii} \not= t_{jj}$, can be rewritten as
\begin{equation} \label{eq:lij}
 \ell_{ij}=-\frac{1}{t_{ii}-t_{jj}}\Big(e_{ij}+\sum_{k=i+1}^n t_{ik}\ell_{ki}-\sum_{k=1}^{j-1}\ell_{i,k}t_{kj}\Big).
\end{equation}
Note that all entries of $L$ appearing on the right-hand side are located either to the left or below the entry $(i,j)$ in the matrix $L$. Thus, if $(i_1,j_1), (i_2,j_2), \ldots, (i_N,j_N)$ with $N = n(n-1)/2$ describes an order in which the entries of $L$
shall be determined, this order must satisfy
\begin{equation} \label{eq:ordercondition}
 \nu < \mu \quad \Longrightarrow \quad i_\nu > i_\mu \text{ or } j_\nu < j_\mu
\end{equation}
for all $\nu,\mu \in \{1,\ldots,N\}$. 
Interestingly, this corresponds to the elimination order of the northeast directed sweep sequences in the nonsymmetric Jacobi algorithm for which
local quadratic convergence was proven in~\cite{Meh08}. 
Each order satisfying~(\ref{eq:ordercondition}) gives rise to a different algorithm for solving~(\ref{eq:linearsystem}). A possible choice is a bottom-to-top columnwise order, leading to Algorithm~\ref{alg:scalar}. Note that we use Matlab notation to refer to entries and submatrices of a matrix.
\begin{algorithm}{Successive substitution for solution of triangular
matrix equation~\eqref{eq:linearsystemagain}}\label{alg:scalar}
 \REQUIRE Upper triangular matrix $T \in\C^{n\times n}$ with pairwise distinct diagonal entries. Strictly lower triangular matrix $E \in\C^{n\times n}$.
 \ENSURE Strictly lower triangular matrix $L \in\C^{n\times n}$ satisfying~(\ref{eq:linearsystemagain}).
 \FOR{$j=1:n-1$}
   \FOR{$i=n:-1:j+1$}
     \STATE $E(i,j) \gets E(i,j)+T(i,i+1:n) \cdot L(i+1:n,j) - L(i,1:j-1) \cdot T(1:j-1,j)$
     \STATE $L(i,j) \gets  -E(i,j)/(T(i,i)-T(j,j))$
   \ENDFOR
 \ENDFOR
\end{algorithm}

\begin{theorem} \label{thm:solvability}
 There exists a unique strictly lower triangular matrix $L$ satisfying~\eqref{eq:linearsystemagain} if and only if the diagonal entries of $T$ are pairwise distinct.
\end{theorem}
\proof
This result can be found (implicitly) in~\cite[Sec. 3]{KonPC94}, \cite{Sun92}, and~\cite[Theorem 3.1]{Sun95schur}. For completeness, we provide an explicit proof.

By construction, Algorithm~\ref{alg:scalar} shows that~\eqref{eq:linearsystemagain} has a solution for every $E$ under stated condition on $T$. The unique solvability simply follows from the fact that~\eqref{eq:linearsystemagain} can be regarded as a linear system of $n(n-1)/2$ equations in 
$n(n-1)/2$ unknowns.

For the other direction, assume that the condition is violated, that is,
$T$ has two identical eigenvalues $\lambda$. Then there is a principal submatrix $\tilde T$ of $T$ taking the form
\[
 \tilde T = \left[ \begin{array}{ccc}
                    \lambda & T_{12} & T_{13} \\
			0 & T_{22} & T_{23} \\
			0 & 0  &\lambda
                   \end{array}\right],
\]
where the eigenvalues of $T_{22}$ are pairwise distinct and different from $\lambda$ or $T_{22}$ vanishes.
We consider 
\[
 \tilde L = \left[ \begin{array}{ccc}
                    0 & 0 & 0 \\
			L_{21} & L_{22} & 0 \\
			1 & L_{32}  &0
                   \end{array}\right],
\]
partitioned conformally with $\tilde T$.
Then $\stril(\tilde T \tilde L - \tilde L \tilde T) = 0$ if and only if
\begin{eqnarray}
 (T_{22} - \lambda I) L_{21} &= & -T_{23} \label{eq:null1} \\
  L_{32} (\lambda I -  T_{22}) &=&  T_{12} \label{eq:null2}\\
 \stril( T_{22}  L_{22}  -  L_{22}  T_{22}) &=& \stril(- T_{23}  L_{32} +  L_{21}  T_{12}) \label{eq:null3}
\end{eqnarray}
Since $\lambda \not\in \Lambda(T_{22})$, equations~(\ref{eq:null1}) and (\ref{eq:null2}) have unique solutions
$L_{21}$ and $L_{32}$. This determines the right-hand side of~(\ref{eq:null3}), which is an equation of the same type as~(\ref{eq:linearsystemagain}). Since $T_{22}$ has pairwise distinct eigenvalues, the first part of this theorem implies that there is a strictly lower triangular solution $L_{22}$ to~(\ref{eq:null3}). By embedding
$L = \text{diag}(0,\tilde L,0)$ we have thus found a nonzero $L \in \stril(\C^{n\times n})$ such that
$\stril(TL - LT) = 0$. Because of linearity, this implies that the equation $\stril(TL - LT) = -E$ is not uniquely solvable for any $E \in \stril(\C^{n\times n})$ if the condition is violated. 
\eproof

The non-locality of the memory access pattern renders an actual implementation of Algorithm~\ref{alg:scalar} slow for larger matrices. Other orderings satisfying~\eqref{eq:ordercondition}, like parallel wavefront techniques~\cite{OleS85}, may lead to more efficient implementations. In the following we use a recursive formulation, as proposed for a variety of matrix equations in~\cite{JonK02a,JonK02b}, to achieve increased data locality in a convenient manner.
For this purpose, let us partition
\begin{equation}\nonumber
  T=\begin{bmatrix}T_{11}&T_{12}\\ 0&T_{22}\end{bmatrix},\quad
  E=\begin{bmatrix}E_{11}&0\\ E_{21}&E_{22}\end{bmatrix},\quad
  L=\begin{bmatrix}L_{11}&0\\ L_{21}&L_{22}\end{bmatrix},
\end{equation}
with $T_{11},E_{11},L_{11}\in\C^{n_1 \times n_1}$. Inserted into~\eqref{eq:linearsystemagain}, this gives
\begin{equation}\label{eq:block2eqn}
\begin{array}{lcl}
 && \begin{bmatrix}\stril(T_{11}L_{11}+T_{12}L_{21}-L_{11}T_{11})&0\\ 
 T_{22}L_{21}-L_{21}T_{11}& \stril(T_{22}L_{22}-L_{21}T_{12}-L_{22}T_{22})\end{bmatrix} \\
 &=&-\begin{bmatrix}E_{11}&0\\ E_{21}&E_{22}\end{bmatrix}.
\end{array}
\end{equation}
The $(2,1)$ entry is a triangular Sylvester equation,
\begin{equation}\label{eq:l21}
  T_{22}L_{21}-L_{21}T_{11}=-E_{21},
\end{equation}
which has a unique solution $L_{21}$, provided that $T_{11}$ and $T_{22}$ have no eigenvalue in common, and can be addressed using, e.g., the software package RECSY~\cite{JonK02a,JonK02b}. Once $L_{21}$ is determined, $L_{11}$ and $L_{22}$ can be obtained as solutions of
\begin{align*}
  \stril(T_{11}L_{11}-L_{11}T_{11}) &=-(E_{11}+\stril(T_{12}L_{21})),\\
  \stril(T_{22}L_{22}-L_{22}T_{22}) &=-(E_{22}-\stril(L_{21}T_{12})).
\end{align*}
Note that both equations are of the same type as~(\ref{eq:linearsystem}), but of size $n_1$ and $n-n_1$, respectively. Their solutions can be obtained by again subdividing into $2\times 2$ blocks or, for smaller problems, by using Algorithm~\ref{alg:scalar}.
The described strategy leads to the following recursive algorithm for solving~\eqref{eq:linearsystem}.

\begin{algorithm}{Recursive block algorithm for triangular matrix equation~\eqref{eq:linearsystemagain}} \label{alg:block} %
\REQUIRE Upper triangular matrix $T \in\C^{n\times n}$ with pairwise distinct diagonal entries. Strictly lower triangular matrix $E \in\C^{n\times n}$. Integer $n_{\min} \ge 2$.
 \ENSURE Strictly lower triangular matrix $L \in\C^{n\times n}$ satisfying~(\ref{eq:linearsystemagain}).
 \IF{$n \le n_{\min}$}
 \STATE Apply Algorithm~\ref{alg:scalar}.
 \ELSE 
 \STATE Choose $n_1$ with $1< n_1 < n$ (e.g., $n_1=\lfloor n/2\rfloor$) and define $\mathsf{ind}_1 = 1:n_1$, $\mathsf{ind}_2 = n_1+1:n$.
 \STATE Solve Sylvester equation \[T(\mathsf{ind}_2,\mathsf{ind}_2)X - XT(\mathsf{ind}_1,\mathsf{ind}_1) = -E(\mathsf{ind}_2,\mathsf{ind}_1).\]
 \STATE Set $L(\mathsf{ind}_2,\mathsf{ind}_1) \gets X$.
 \STATE \label{eq:mmult1} Update $E(\mathsf{ind}_1, \mathsf{ind}_1) \gets E(\mathsf{ind}_1, \mathsf{ind}_1) + {\rm stril}(  T(\mathsf{ind}_1,\mathsf{ind}_2) L(\mathsf{ind}_2,\mathsf{ind}_1)   )$. 
 \STATE \label{eq:mmult2} Update $E(\mathsf{ind}_2, \mathsf{ind}_2) \gets E(\mathsf{ind}_2, \mathsf{ind}_2) - {\rm stril}(  L(\mathsf{ind}_2,\mathsf{ind}_1) T(\mathsf{ind}_1,\mathsf{ind}_2)   )$. 
 \STATE Use Algorithm~\ref{alg:block} with input $T(\mathsf{ind}_1,\mathsf{ind}_1)$, $E(\mathsf{ind}_1,\mathsf{ind}_1)$ to compute $L(\mathsf{ind}_1,\mathsf{ind}_1)$.
 \STATE Use Algorithm~\ref{alg:block} with input $T(\mathsf{ind}_2,\mathsf{ind}_2)$, $E(\mathsf{ind}_2,\mathsf{ind}_2)$ to compute $L(\mathsf{ind}_2,\mathsf{ind}_2)$.
\label{line:rhsupdate}
 
 \ENDIF
\end{algorithm}

The main cost of Algorithm~\ref{alg:block} is in the two matrix multiplications in Lines~\ref{eq:mmult1} and~\ref{eq:mmult2}, which allows to leverage the efficiency of level 3 BLAS operations. 

\subsection{Orthogonalization procedures} \label{sec:orthog}

In the following, we discuss algorithms for carrying out the (approximate) orthogonalization in lines~\ref{line:optqr} and~\ref{line:unitary} of Algorithm~\ref{alg:template}.
A suitable orthogonalization procedure needs to satisfy two goals. On the one hand, it should improve orthogonality. On the other hand, it should not modify the matrix more than needed. With these two goals in mind, we require the following: For $\varepsilon > 0$, consider any matrices $Q,W \in \C^{n\times n}$ such that $W$ is skew-Hermitian and 
\begin{equation} \label{eq:orthentry}
 \|Q^H Q - I\|_F \le c_Q \varepsilon^2, \quad \|W\|_F \le c_W \varepsilon,
\end{equation}
with constants $c_Q, c_W$ independent of $\varepsilon$. Then the matrix $Q_{\mathsf{new}}$ returned by the orthogonalization procedure applied to $Q(I+W)$ satisfies
\begin{equation} \label{eq:orthexit}
 \| Q_{\mathsf{new}}^H Q_{\mathsf{new}} - I\|_F = O(\varepsilon^4), \quad \|Q_{\mathsf{new}} - Q(I+W) \|_F = O(\varepsilon^2).
\end{equation}

\subsubsection{Orthogonalization by QR decomposition}

When $Q$ is unitary (that is,~\eqref{eq:orthentry} is satisfied with $c_Q = 0$) then $Q(I+W)$ is an element of the tangent space of the manifold of unitary matrices at $Q$. Retractions, a concept popularized in the context of Riemannian optimization~\cite{Absil2008,boumal2022intromanifolds}, map an element from the tangent space back to the manifold. Thus, the first relation in~\eqref{eq:orthexit} is trivially satisfied by a retraction. The second relation in~\eqref{eq:orthexit} follows from the fact that retractions approximate the exponential map in first order.
Suitable retractions and, hence, suitable orthogonalization procedures are obtained from extracting the orthogonal factors of a QR or polar decomposition~\cite[Sec. 7.3]{boumal2022intromanifolds}.  More specifically, if one computes a QR decomposition
\[
  I+W = Q_L R,
\]
such that $Q_L$ is unitary and the upper triangular matrix $R$ has real and positive entries on its diagonal then $Q_{\mathsf{new}} = Q Q_L$ satisfies~\eqref{eq:orthexit}.

\subsubsection{Orthogonalization by Newton-Schulz iteration}

The Newton-Schulz iteration $Q_{k+1} = \frac12 Q_k (3I - Q_k^H Q_k)$ converges quadratically towards the unitary factor of the polar decomposition of $Q_0$ if $Q_0$ is invertible and $\|Q_0\|_2 < \sqrt{3}$~\cite[Ch. 8]{Hig08}. The following lemma shows that one step of the Newton-Schulz iteration yields a suitable orthogonalization procedure.
\begin{lemma} 
\label{lemma::NS}
Consider $\hat Q = Q (I+W)$ satisfying~\eqref{eq:orthentry}. Then
\[
 Q_{\mathsf{new}} = \frac12 \hat Q (3I - \hat Q^H \hat Q)
\]
satisfies~\eqref{eq:orthexit}.
\end{lemma}
\proof
The first relation of~\eqref{eq:orthexit} follows from 
$Q_{\mathsf{new}}^H Q_{\mathsf{new}} - I = -\frac34 (\hat Q^H \hat Q - I)^2 + \frac14 (\hat Q^H \hat Q - I)^3$; see~\cite[Problem 8.20]{Hig08}.
To show the second relation, we set $\triangle := Q^H Q - I$ and compute
\[
 \hat Q^H \hat Q - I = (I-W)(I+\triangle)(I+W) - I = \triangle -W^2 -W \triangle + \triangle W - W\triangle W.
\]
Combined with $Q_{\mathsf{new}} = \hat Q -\frac12 \hat Q (\hat Q^H \hat Q -I)$,
this shows
\[
 \|Q_{\mathsf{new}} - \hat Q\|_F \le \frac{c_Q + c_W^2}{2} \varepsilon^2 + O(\varepsilon^3).
\]
\eproof

\subsection{Local convergence analysis}

We now perform a local convergence analysis of Algorithm~\ref{alg:template}. This analysis makes use of the quantity
\begin{equation} \label{eq:phi}
\phi(T):=\min \big\{ \|\stril(TL-LT)\|_F: L\in \stril(\C^{n\times n}), \|L\|_F=1 \big\}
\end{equation}
for an upper triangular matrix $T$, where $\stril(\C^{n\times n})$ denotes the set of all strictly lower triangular $n\times n$ matrices. Note that $1/\phi(T)$ governs the first-order sensitivity of a Schur decomposition~\cite{KonPC94,Sun92}.
By Theorem~\ref{thm:solvability}, $\phi(T) > 0$ if and only if the diagonal entries of $T$ are pairwise distinct. 
The following lemma shows that $\phi$ is Lipschitz continuous with constant $\sqrt{2}$; see~\cite[Theorem 3.2]{Sun95schur} for a related result.

\begin{lemma} \label{lemma:lipschitzphi}
Consider upper triangular matrices $T,\hat T$, each having pairwise distinct diagonal entries.
Then $|\phi(T)-\phi(\hat{T})|\le 2 \|T-\hat T\|_2$.
\end{lemma}
\proof
By the definition~\eqref{eq:phi}, the quantity $\phi(T)$ is the smallest singular value of the linear operator $\mathcal L_T: \stril(\C^{n\times n})\to \stril(\C^{n\times n})$, $\mathcal L_T(L): L \mapsto \stril(TL-LT)$, with the norm $\|\mathcal L_T\| = \sup\{ \|\mathcal L_T(L)\|_F: \|L\|_F = 1, L \in \stril(\C^{n\times n})\}$. From 
\begin{eqnarray*}
  \|(\mathcal L_T - \mathcal L_{\hat T})(L)\|_F &=& \|\stril((T-\hat T)L-L(T-\hat T))\|_F \\ & \le & 
  \|(T-\hat T)L-L(T-\hat T)\|_F 
  \le 2 \|T-\hat T\|_2 \|L\|_F
\end{eqnarray*}
it follows that $\|\mathcal L_T - \mathcal L_{\hat T}\| \le 2 \|T-\hat T\|_2$. Since Weyl's inequalities~\cite[Corollary 7.3.5]{Horn2013} imply that singular values are Lipschitz continuous with constant $1$, this concludes the proof.
\eproof

Suppose that $A$ has mutually distinct eigenvalues. Then its Schur decomposition $A = QTQ^H$ is unique up to the order of the eigenvalues on the diagonal of $T$ and unitary diagonal transformations. While $\phi(T)$ is invariant under the latter, it does depend on the eigenvalue order.
To circumvent this difficulty, we introduce
\begin{equation} \label{eq:phia}
\psi(A):=\min\{ \phi(T): \text{$A = QTQ^H$ is a Schur decomposition}\}. 
\end{equation}
Because there are only finitely many different eigenvalue orderings, the minimum is assumed and positive. The following lemma shows that $\psi(\cdot)$ remains positive in a neighborhood of $A$.

\begin{lemma} \label{eq:s}
Let $A \in \C^{n\times n}$ have mutually distinct eigenvalues and consider $\varepsilon > 0$ sufficiently small. Then for all $\tilde A\in \C^{n\times n}$ with $\|\tilde A - A\|_F \le \varepsilon$, it holds that $\psi(\tilde A) \ge
\psi(A) - 2 (1+ 4 \|A\|_2 / \psi(A))\varepsilon$.
\end{lemma}
\proof
For fixed $\triangle \in \C^{n\times n}$ with $\|\triangle\|_F = 1$ and $\varepsilon > 0$, set $A(\varepsilon):= A + \varepsilon\triangle$ and consider a Schur decomposition $A(\varepsilon) = Q(\varepsilon) T(\varepsilon) Q(\varepsilon)^H$. Let $\lambda_1(\varepsilon),\ldots,\lambda_n(\varepsilon)$ denote the diagonal elements of $T(\varepsilon)$. Because of continuity of eigenvalues, we may assume w.l.o.g. that each $\lambda_i(\cdot)$ is a continuous function.

For sufficiently small $\varepsilon > 0$, the perturbation results from~\cite{KonPC94,Sun92,Sun95schur} imply that there is a Schur decomposition $A(\varepsilon) = \tilde Q(\varepsilon) \tilde T(\varepsilon) \tilde Q(\varepsilon)^H$ with $\tilde T(\varepsilon) \approx T(0)$. Specifically, the computation of the condition number of $T$ on Page 391 of~\cite{KonPC94} shows that
\[
\|\tilde T(\varepsilon) - T(0)\|_F \le \Big(1+2 \sqrt{2} \frac{\|T(0)\|_2}{\phi(T(0))} \Big) \varepsilon + O(\varepsilon^2) \le 
\Big(1+ 4 \frac{\|T(0)\|_2}{\phi(T(0))} \Big) \varepsilon.
\]
By Lemma~\ref{lemma:lipschitzphi}, it follows that
\begin{equation} \label{eq:lowerbound}
 \phi(\tilde T(\varepsilon)) \ge \phi(T(0)) - 2 \Big(1+ 4 \frac{\|T(0)\|_2}{\phi(T(0))} \Big) \varepsilon \ge
 \psi(A) - 2 \Big(1+ 4 \frac{\|A\|_2}{\psi(A)} \Big) \varepsilon.
 \end{equation}
The perturbation bound on $\tilde T(\varepsilon)$ also implies that the eigenvalues $\lambda_1(\varepsilon),\ldots,\lambda_n(\varepsilon)$ need to appear in this order on the diagonal of $\tilde T(\varepsilon)$ for sufficiently small $\varepsilon>0$. Thus, 
$\tilde T(\varepsilon)$ and $T(\varepsilon)$ only differ by a transformation with a unitary diagonal matrix and therefore
$\phi(\tilde T(\varepsilon)) = \phi(T(\varepsilon))$. Combined with~\eqref{eq:lowerbound}, this completes the proof. 
\eproof

\begin{theorem} \label{thm:mainconvergence}
Assume that the eigenvalues of $A\in \C^{n\times n}$ are pairwise distinct.
Given $\varepsilon > 0$ consider any matrix $Q\in \C^{n\times n}$ satisfying \begin{enumerate}
 \item $\|\stril(Q^H A Q)\|_F \le \varepsilon$;
\item $\|Q^H Q - I\|_F \le \varepsilon^2$.
\end{enumerate}
Let $Q_{\mathsf{new}}$ denote the output of one iteration of Algorithm~\ref{alg:template} applied to $A,Q$ using an orthogonalization procedure that satisfies~\eqref{eq:orthexit}. Then, for $\varepsilon > 0$ sufficiently small, it holds that
\begin{equation} \label{eq:claims}
 \|\stril(Q_{\mathsf{new}}^H A Q_{\mathsf{new}})\|_F = O(\varepsilon^2), \quad \|Q_{\mathsf{new}}^H Q_{\mathsf{new}} - I\|_F = O(\varepsilon^4).  \quad 
\end{equation}
\end{theorem}
\proof
Let $E = \stril(Q^H A Q)$ and $T = Q^H A Q -E$. By the definition of $\phi(T)$, it follows that 
the solution $L$ of~(\ref{eq:linearsystem}) satisfies
\[
 \|L\|_F\le \frac{\|E\|_F}{\phi(T)} \le \frac{\varepsilon}{\phi(T)}.
\]
In order to proceed from here, we need to show that $\phi(T)$, which depends on $Q$, admits a uniform lower bound.
By the polar decomposition there is a unitary matrix $\tilde Q$ such that
$
 \|\tilde Q - Q\|_F \le \|Q^H Q - I\|_F \le \varepsilon^2;
$
see, e.g.,~\cite[P. 380]{Hig02}.
In turn, we obtain the perturbed Schur decomposition
\begin{eqnarray*}
 \tilde Q^H A \tilde Q &=& Q^H A Q + \tilde Q^H A(\tilde Q-Q) + (\tilde Q - Q)^H A Q \\
 &=& T+ E + \triangle, \quad \|\triangle\|_F \le 3 \|A\|_2\varepsilon^2,
\end{eqnarray*}
where we used that $\|Q\|_2 = \sqrt{\|Q^H Q\|_2} \le \sqrt{1 +  \varepsilon^2} \le 2$ for $\varepsilon$ sufficiently small. Equivalently, this can be viewed as a Schur decomposition of a perturbed matrix: $A - \tilde Q(E+\triangle) \tilde Q^H = \tilde Q T \tilde Q^H$. From Lemma~\ref{eq:s} it follows that
$
 \phi(T) \ge \psi(A)/2 
$
holds for $\varepsilon$ sufficiently small. 
Hence, $W = L-L^H$ satisfies $\|W\|_F \le c_W \varepsilon$ with $c_W = 4 / \psi(A)$. This means that the matrix $\hat Q = Q(I+W)$ satisfies the conditions~\eqref{eq:orthentry} and therefore the matrix $Q_{\mathsf{new}}$ returned by the orthogonalization procedure in line~\ref{line:unitary} of Algorithm~\ref{alg:template} satisfies~\eqref{eq:orthexit}. This already establishes the second relation in~\eqref{eq:claims}. To establish the first relation, we note that 
\begin{eqnarray*}
 \hat Q^H A \hat Q &=& (I + W)^H ( T + E ) (I + W) \\
&=& \underbrace{E - L T + T L}_{\text{upper triangular}} + \underbrace{T + L^H T-TL^H}_{\text{upper triangular}} + W^H(T+E) W + W^H E + EW,
\end{eqnarray*}
and, hence,
\begin{eqnarray*}
  \big\|\stril(\hat Q^H A \hat Q)\big\|_F &\le& \|T+E\|_F \|W\|_F^2 +2\|W\|_F\|E\|_F  \\
  &\le & \|Q\|_2^2 \|A\|_F c_W^2 \varepsilon^2 + 2c_w\varepsilon^2 \le (4 \|A\|_F c_W^2 + 2c_W) \varepsilon^2.
\end{eqnarray*}
Combined with $\|Q_{\mathsf{new}}-\hat Q\| = O(\varepsilon^2)$, this proves the first relation in~\eqref{eq:claims}.
\eproof

If Algorithm~\ref{alg:template} is started with a unitary matrix $\hat Q$ sufficiently close to a unitary matrix $Q$ that transforms $A$ to Schur form then the relations~\eqref{eq:orthexit} imply that the matrix $Q$ obtained from the orthogonalization procedure applied in line~\ref{line:optqr} satisfies the conditions of Theorem~\ref{thm:mainconvergence}. Hence, Theorem~\ref{thm:mainconvergence} establishes local quadratic convergence of Algorithm~\ref{alg:template}.

\begin{remark} \label{remark:landing} The idea of imposing constraints (such as unitary matrix structure) only asymptotically upon convergence of an iterative procedure is not new. Similar ideas have been proposed in the context of differential-algebraic equations (Baumgarte's method~\cite{Ascher1994}), constrained optimization (interior point method~\cite{Nocedal2006}), and -- more recently -- Riemannian optimization (\cite{Gao2018}, landing algorithm~\cite{Ablin2021}). However, Algorithm~\ref{alg:template} does not seem to fit into any of these existing developments.
\end{remark}

\subsection{Complete algorithm for mixed precision} \label{sec:completealg}

In this section we specialize the template Algorithm \ref{alg:template} to computing a Schur decomposition of a given matrix $A$ in a mixed precision environment. The decomposition is first computed in a lower precision (lp in the following), and then refined to a higher precision (hp in the following). A typical scenario uses double precision as lp, and quadruple or $100$-digits precision as hp.
The initial lp decomposition produces the matrices $\hat{Q}$ and $T_\textsf{lp}$ such that $A \approx \hat Q T_{\textsf{lp}} \hat Q^H$, and the matrices $A$ and $\hat Q$ are the input to the refinement algorithm. In order to ensure convergence (see Theorem~\ref{thm:mainconvergence} and the discussion in the previous section), orthogonality of the matrix $\hat Q$ is first improved by applying one step of the Newton-Schulz iteration:
\begin{equation}
    \label{alg:eq:NS}
    Q = \frac{1}{2} \hat{Q}(3I - \hat{Q}^H \hat{Q}).
\end{equation}
This computation is done entirely in hp; the matrix $\hat{Q}$ is first converted to hp. The most time consuming parts of \eqref{alg:eq:NS} are the two matrix-matrix multiplications in hp.

We then enter the loop in Algorithm~\ref{alg:template}: computation of $\hat{T} = Q^H A Q$ in Line~\ref{line:extractT} is also done in hp, requiring two additional matrix-matrix multiplications. Equation~\eqref{eq:linearsystem} in Line~\ref{line:solvebasic} can be solved entirely in lp. Note that numerical instabilities in computing $L$ are expected when there are clustered eigenvalues in $T$. In hope that these instabilities do not spread through the entire matrix $L$, the initial lp Schur decomposition is reordered so that clustered eigenvalues appear on neighboring diagonal elements of $T_{\textsf{lp}}$. This is done using \code{ordschur} in Matlab; the target permutation of eigenvalues is determined by the order in which they appear after being projected on a random line in the complex plane.

Finally, we need to compute the matrix $Q_{\textsf{new}}$ by improving the orthogonality of $Q(I+L-L^H)$.
In line with Lemma~\ref{lemma::NS}, this is done by applying one step of the Newton-Schulz iteration to $Q(I+L-L^H)$. To lower computational complexity and  avoid unnecessary matrix-matrix multiplications in hp, we proceed in the following way:
let $W = L-L^H$, and $Y = Q^H Q - I$. Then \eqref{alg:eq:NS} with $\hat{Q} = Q(I+W)$ becomes
\begin{align}
    Q_{\textsf{new}}
        &= \frac{1}{2} Q(I+W) (3I - (I+W)^H Q^H Q (I+W)) \nonumber \\
        &= \frac{1}{2}Q(2I + 2W - Y - YW + W^2 + W^3 + W^2Y + W^2 Y W), \label{alg:eq:NSmerged}
\end{align}
where we used $W^H = -W$. We observe that the products between the matrices $Y$ and $W$ can all be computed in lp while the summation of the terms needs to be carried out in hp.
In any case, only $2$ hp matrix-matrix multiplications are needed to compute the update: one to compute $Y$, and one to compute the product of $Q$ with $(2I + 2W - Y - \cdots)$. The cost can be reduced slightly by observing that high-order terms tend become very small and may even fall below the unit roundoff in lp. In our implementation, we actually use the approximation $Q_{\textsf{new}} \approx \frac{1}{2}Q(2I + 2W - Y - YW + W^2 + W^3)$.

The complete procedure is summarized in Algorithm~\ref{alg:mixedprec}. As the numerical experiments will demonstrate, the computational complexity is almost entirely concentrated in hp matrix multiplication. The algorithm does $2$ such multiplications before the loop, and then $4$ in each loop iteration. The iteration is stopped in line~\ref{line:extractE} when $E$ becomes small enough. Note that lines~\ref{line:sorth}--\ref{line:eorth} can be skipped if the norms of $Y$ and $W$ indicate that $Q(I+W)$ is already unitary in hp. For this purpose, the norm of $\|Y\|_F \sim \mathsf{u}_{\mathsf{hp}}$,  where $\mathsf{u}_{\mathsf{hp}}$ denotes unit roundoff in hp, while $\|W\|_F \sim \sqrt{\mathsf{u}_{\mathsf{hp}}}$ due the fact that $W$ is skew-Hermitian (see~\eqref{eq:orthW}). This potentially saves 1 hp matrix multiplication in the penultimate iteration.

\begin{algorithm}{Mixed-precision computation of Schur decomposition}\label{alg:mixedprec}
\REQUIRE Matrix $A \in \C^{n\times n}$ in hp.
\ENSURE Matrix $Q \in \C^{n\times n}$ such that $Q$ is unitary in hp and $Q^H A Q$ is upper triangular in hp.
 \STATE Compute (ordered) Schur decomposition $A \approx \hat Q T_{\textsf{lp}} \hat Q^H$ in lp.
 \STATE Convert $\hat Q$ to hp and update $Q \gets \frac{1}{2} \hat{Q}(3I - \hat{Q}^H \hat{Q})$ in hp (one step of Newton-Schulz).
 \REPEAT
   \STATE Compute $\hat T \gets Q^H A Q$ in hp.
   \STATE \label{line:extractE} Set $E\leftarrow \stril(\hat T)$, $T\leftarrow \hat T - E$ and convert to lp.
   \STATE \label{line:solve}Solve equation~(\ref{eq:linearsystem}) for strictly lower triangular $L$ in lp, using Algorithm~\ref{alg:block}.
   \STATE Convert $L$ to hp.
   \STATE Set $W = L-L^H$, compute $Y = Q^H Q - I$ in hp.
   \STATE \label{line:sorth}  Compute products $YW$, $W^2$, $W^3$, $W^2 Y$, $W^2 YW$ in lp.
   \STATE Compute $\Sigma = 2I + 2W - Y - YW + W^2 + W^3$ in hp. 
   \STATE \label{line:eorth} Update $Q \gets \frac{1}{2}Q \Sigma$ in hp (one step of Newton-Schulz).
\UNTIL{convergence}
\end{algorithm}

\subsection{Symmetric $A$} \label{sec:symmetrica}

When $A$ is real and symmetric, its Schur form becomes diagonal and our algorithms can be simplified significantly. Instead of Steps 4 of Algorithms~\ref{alg:template} and~\ref{alg:mixedprec}, $\hat T$ is now decomposed in its diagonal part $T$ and its off-diagonal part $E$. In turn, the solution $L$ of the linear system~\eqref{eq:linearsystem} becomes nearly trivial; its entries are given by $\ell_{ij} = e_{ij} / ( t_{ii} - t_{jj})$ for $i > j$. The resulting simplified algorithms bear close similarity with the method {\sf RefSyEv} proposed by Ogita and Aishima~\cite{Ogita18}. However, in contrast to our approach, {\sf RefSyEv} performs the improvement of diagonality and orthogonality in reverse order and integrates them in a single update. More concretely, using our notation, one iteration of {\sf RefSyEv} updates $Q \gets Q(I+W)$ in the following manner.
It first computes $Y = Q^T Q - I$ and $\hat T = Q^T A Q$. Orthogonality is improved by setting the symmetric part of $W$ to $-Y/2$, which is equivalent to one step of Newton-Schulz as pointed out in~\cite[Appendix A]{Ogita18}. Then the diagonal of $\hat T$ is updated and the skew-symmetric part of $W$ is chosen to improve diagonality by solving an equation of the form~\eqref{eq:linearsystem}, with the right-hand side updated to reflect improvement of orthogonality. {\sf RefSyEv} actually does not treat the symmetric and skew-symmetric parts separately but derives a simple and elegant formula to update all entries of $W$ in a single step.  The most significant advantage of {\sf RefSyEv} is that it avoids the initial orthogonalization in Algorithms~\ref{alg:template} and~\ref{alg:mixedprec}. This saves $2$ hp matrix-matrix multiplications in the initial phase; note, however,  that both {\sf RefSyEv} and Algorithm~\ref{alg:mixedprec} need the same number ($4$) of hp matrix-matrix multiplications in each iteration. In summary, there seems to be little that speaks in favor of our approach for symmetric eigenvalue problems, especially when taking into account that 
{\sf RefSyEv} and its further development described in~\cite{Ogita2019} take precautions for (nearly) multiple eigenvalues in order to still attain fast convergence in such a critical case. One (small) advantage of our approach is that the separation of the orthogonalization step allows for the use of other orthogonalization procedures, which may be of interest for future developments.


\section{Numerical experiments}

In this section, we report the results of a number of numerical experiments in order to demonstrate the correctness and effectiveness of the proposed algorithm. The experiments were executed on a notebook computer with Intel Core i5-1135G7 CPU and 24GB RAM, running Ubuntu 21.10. Algorithm~\ref{alg:mixedprec} was implemented in Matlab R2021b. We have also implemented a straightforward extension of this algorithm to real Schur decompositions~\cite{Golub2013}. The main difference to the complex case is that some attention needs to be paid to $2\times 2$ diagonal blocks containing complex pairs of eigenvalues. In Algorithm~\ref{alg:block}, the value of $n_{\text{min}}$ was set to $4$ in the complex case, and to $1$ in the real case. To solve Sylester equations of the form~\eqref{eq:l21} we used the internal Matlab solver \code{matlab.internal.math.sylvester\_tri}.

Our test procedure starts by generating an input matrix $A$ in high precision, for which we discuss two scenarios: quadruple precision (i.e., 34 decimal digits of precision), and 100 decimal digits of precision. The matrix is given as input to Algorithm~\ref{alg:mixedprec}. There we use the standard double precision as low precision, and the builtin Matlab command \code{schur} to compute the initial lp Schur decomposition (either real or complex).
The refined factors $Q$, $T$ returned by Algorithm~\ref{alg:mixedprec} are verified for correctness: we check if $T$ has proper (quasi) triangular form, analyze orthogonality of $Q$ by computing $\|I - Q^H Q\|_F$, and compute the residual $\|T - Q^H A Q\|_F$.

The performance of our algorithm heavily relies on the efficiency of matrix-matrix multiplication in the target precision. For that purpose we use the Matlab toolbox \code{acc} based on the Ozaki scheme~\cite{Ozaki12}, which was also used in \cite{Ogita18}. The toolbox uses a configurable number of standard doubles to represent a single high-precision number; we use 2 doubles to represent a quadruple precision number, and 6 doubles to represent a 100-digit number.

The time needed by the whole refinement procedure (including the initial Schur decomposition in double precision) is compared with the time needed for computing the Schur decomposition immediately in quadruple/100 digits precision. For the latter we use Advanpix \cite{advanpix}, a multiprecision computing toolbox for Matlab.\footnote{Note that
Matlab's Symbolic Toolbox vpa (variable precision arithmetic) supports eigenvalue computations but it does not support the computation of Schur decompositions. Executing {\tt eig} in vpa takes $33.70$s for a $100 \times 100$ matrix in quadruple precision, compared to only $0.30$s needed by Advanpix for the complete Schur decomposition.}

\begin{example}
	This first example aims at verifying the accuracy of our algorithm by applying it to the companion matrix  for the Wilkinson polynomial
	$$
		p_{20}(x) = \prod_{i=1}^{20} (x - i).
	$$ 
	The coefficients of the polynomial and, hence, entries of the companion matrix vary wildly in magnitude and cannot be stored exactly in double precision. 
	As a consequence, computing the eigenvalues via the Matlab command \code{roots(poly(1:20))} incurs a large error; see Table~\ref{tbl:example1-wilkinson}.

	\begin{table}
		\begin{center}
			{\tiny
			\begin{tabular}{|r|r|r|r|r|}
				\hline
				exact & \code{roots(poly(1:20))} & \makecell{Schur in high-precision, \\ Advanpix with 34 digits} & \makecell{Algorithm~\ref{alg:mixedprec} (double $\to$ quad, \\ using Advanpix with 34 digits)} & \makecell{Algorithm~\ref{alg:mixedprec} (double $\to$ quad, \\ using \code{acc} with 2 cells)} \\ \hline
				20    &  $1.25\cdot 10^{-4}$                &  $-2.55\cdot 10^{-23}$    &  $-2.19\cdot 10^{-23}$  &  $-1.02\cdot 10^{-21}$    \\ 
				19    &  $-1.29\cdot 10^{-3}$                &  $8.41\cdot 10^{-22}$    &  $-2.70\cdot 10^{-22}$  &  $9.17\cdot 10^{-21}$    \\
				18    &  $6.32\cdot 10^{-3}$                &  $-6.20\cdot 10^{-21}$    &  $2.25\cdot 10^{-21}$  &  $-3.67\cdot 10^{-20}$    \\
				17    &  $-1.85\cdot 10^{-2}$                &  $2.27\cdot 10^{-20}$    &  $-7.09\cdot 10^{-21}$  &  $9.03\cdot 10^{-20}$    \\
				16    &  $4.02\cdot 10^{-2}$                &  $-5.08\cdot 10^{-20}$    &  $1.31\cdot 10^{-20}$  &  $-1.60\cdot 10^{-19}$    \\
				15    &  $-5.93\cdot 10^{-2}$                &  $7.55\cdot 10^{-20}$    &  $-1.67\cdot 10^{-20}$  &  $2.27\cdot 10^{-19}$    \\
				14    &  $6.98\cdot 10^{-2}$                &  $-7.72\cdot 10^{-20}$    &  $1.58\cdot 10^{-20}$  &  $-2.66\cdot 10^{-19}$    \\
				13    &  $-6.26\cdot 10^{-2}$                &  $5.40\cdot 10^{-20}$    &  $-1.15\cdot 10^{-20}$  &  $2.52\cdot 10^{-19}$    \\
				12    &  $4.11\cdot 10^{-2}$                &  $-2.40\cdot 10^{-20}$    &  $6.47\cdot 10^{-21}$  &  $-1.86\cdot 10^{-19}$    \\
				11    &  $-2.24\cdot 10^{-2}$                &  $4.73\cdot 10^{-21}$    &  $-2.71\cdot 10^{-21}$  &  $1.03\cdot 10^{-19}$    \\
				10    &  $8.80\cdot 10^{-3}$                &  $1.68\cdot 10^{-21}$    &  $8.29\cdot 10^{-22}$  &  $-4.24\cdot 10^{-20}$    \\
				 9    &  $-2.71\cdot 10^{-3}$                &  $-1.70\cdot 10^{-21}$    &  $-1.70\cdot 10^{-22}$  &  $1.25\cdot 10^{-20}$    \\
				 8    &  $6.05\cdot 10^{-4}$                &  $6.72\cdot 10^{-22}$    &  $2.59\cdot 10^{-23}$  &  $-2.62\cdot 10^{-21}$    \\
				 7    &  $-9.69\cdot 10^{-5}$                &  $-1.58\cdot 10^{-22}$    &  $-2.48\cdot 10^{-24}$  &  $3.76\cdot 10^{-22}$    \\
				 6    &  $1.04\cdot 10^{-5}$                &  $2.33\cdot 10^{-23}$    &  $1.60\cdot 10^{-25}$  &  $-3.46\cdot 10^{-23}$    \\
				 5    &  $-7.05\cdot 10^{-7}$                &  $-2.07\cdot 10^{-24}$    &  $-1.18\cdot 10^{-26}$  &  $1.70\cdot 10^{-24}$    \\
				 4    &  $2.61\cdot 10^{-8}$                &  $1.02\cdot 10^{-25}$    &  $9.03\cdot 10^{-28}$  &  $-2.17\cdot 10^{-26}$    \\
				 3    &  $-4.44\cdot 10^{-10}$                &  $-2.44\cdot 10^{-27}$    &  $-6.45\cdot 10^{-29}$  &  $-1.14\cdot 10^{-27}$    \\
				 2    &  $1.61\cdot 10^{-12}$                &  $2.85\cdot 10^{-29}$    &  $-1.30\cdot 10^{-30}$  &  $-8.12\cdot 10^{-28}$    \\
				 1    &  $5.11\cdot 10^{-14}$                &  $-1.95\cdot 10^{-31}$    &  $-2.64\cdot 10^{-29}$  &  $-8.16\cdot 10^{-28}$    \\ \hline

			\end{tabular}}		
		\end{center}
		\caption{First column: Eigenvalues of the companion matrix for the Wilkinson polynomial. Other columns: Absolute errors when eigenvalues are computed in different ways.}
		\label{tbl:example1-wilkinson}
	\end{table}
	Storing the companion matrix in quadruple precision allows for more accurate eigenvalues.
	Table~\ref{tbl:example1-wilkinson} compares the accuracy of the results obtained from applying Advanpix's {\tt schur} (with 34 digits) with the ones from Algorithm~\ref{alg:mixedprec} using mixed double-quadruple precision. In both cases the obtained accuracy is on a similar level and significantly improves upon the double precision computation. However, it can also be seen that the use of the \code{acc} toolbox (which uses 2 doubles representing each number) results in a slight loss of accuracy compared to using Advanpix (34 digits) within Algorithm~\ref{alg:mixedprec}.
	A similar effect is seen when using 100 digits in Advanpix versus 6 doubles per high-precision number in \code{acc}. As the \code{acc} toolbox is significantly faster, we will use it in all subsequent experiments. Note that ongoing and future modifications of the Ozaki scheme, such as the ones presented in~\cite{Lange2022}, will likely yield further improvements of the accuracy and efficiency of this approach.
\end{example}

\begin{example} \label{example:performance}
	This experiment focuses on the performance of Algorithm~\ref{alg:mixedprec}. We generate a series of random matrices of increasing sizes with random entries from the standard normal distribution. Both the real and the complex Schur forms are computed in high precision. As can be seen from Figure \ref{fig:example2-speedup},
	Algorithm~\ref{alg:mixedprec} is much faster than Advanpix's {\tt schur}. For this example, our algorithm always converges within $3$ iterations to quadruple precision, while $7$--$8$ iterations are needed to attain $100$-digit precision. The computed factors $Q$ and $T$ satisfy
	\begin{align*}
		\|I-Q^\ast Q\|_F & \leq \left\{ 
			\begin{array}{ll} 
				9 \cdot 10^{-32},& \text{in quadruple precision;} \\
				3 \cdot 10^{-97},& \text{in $100$-digits precision;} 
			\end{array} \right. \\
		\|\stril(Q^\ast A Q)\|_F / \|A\|_F & \leq \left\{ 
			\begin{array}{ll} 
				3 \cdot 10^{-33},& \text{in quadruple precision;} \\
				2 \cdot 10^{-98},& \text{in $100$-digits precision,} 
			\end{array} \right.
	\end{align*}
	for all input matrices $A$.

	\begin{figure}
		\centering
		\includegraphics[width=0.75\textwidth]{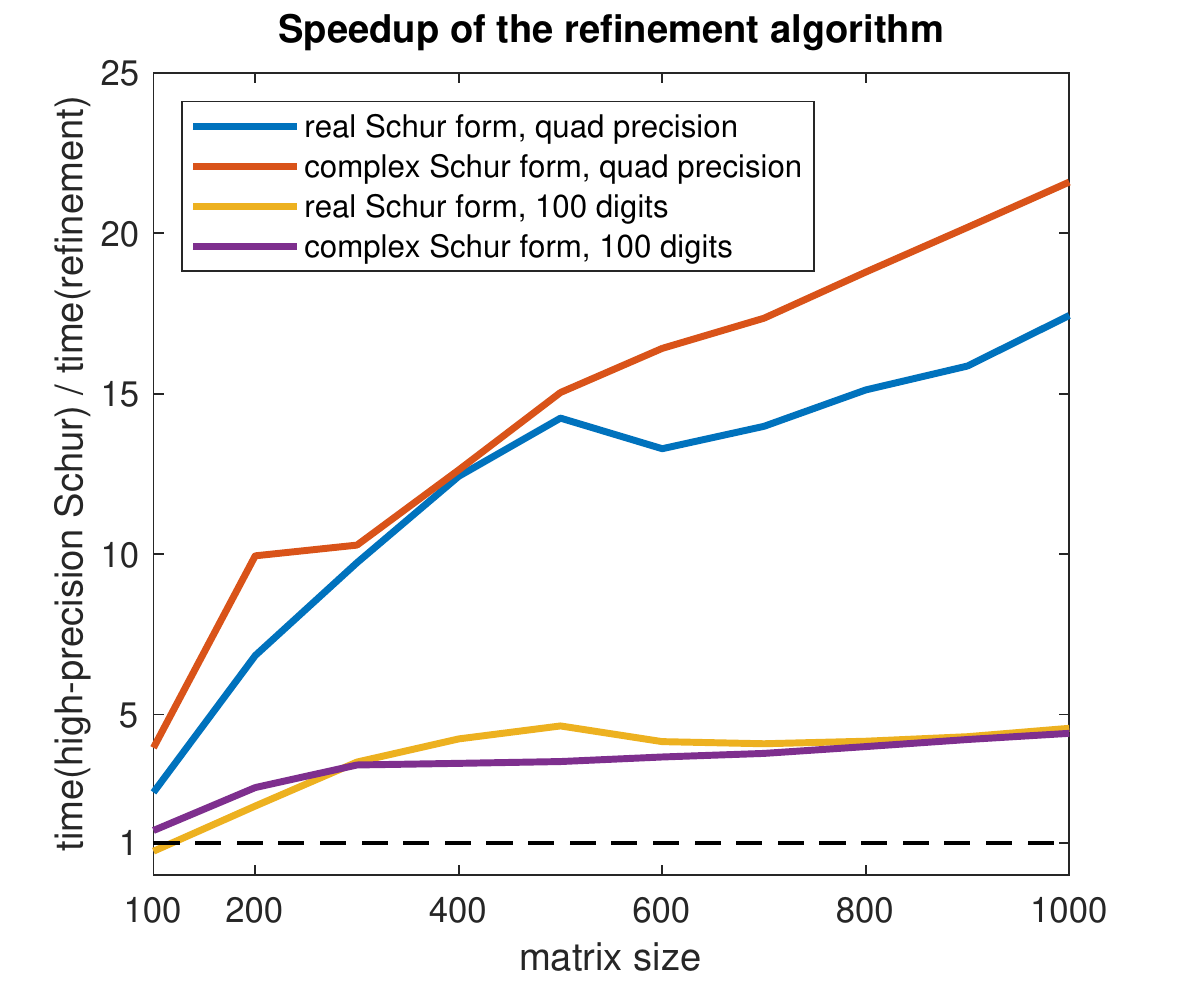}
		\caption{Speedup attained by Algorithm~\ref{alg:mixedprec} for Example~\ref{example:performance}. Time needed by Advanpix
		for computing high-precision Schur decomposition divided by time needed by Algorithm~\ref{alg:mixedprec} (including the time needed to compute Schur decomposition in double precision).}
		\label{fig:example2-speedup}
	\end{figure}

	For more insight, Table \ref{tbl:example2-timings} shows detailed timings for the largest matrix of size $1000$. We note that  Algorithm~\ref{alg:mixedprec} spends essentially all of its time on high-precision matrix-matrix multiplications.

	\begin{table}
		\begin{center}
			{\small
			\begin{tabular}{|c||c|c|c|c|}
				\hline
					                       & quad, real        & quad, complex      & 100 digits, real    & 100 digits, complex \\ \hline
					high-precision Schur   & $165.14$s         & $614.61$s          & $713.66$s           & $2178.66$s \\ \hline
					\makecell{refinement algorithm: \\ total time}             
					                       & $9.47$s           & $28.45$s           & $156.18$s           & $493.72$s \\ \cdashline{2-5}
					matrix multiplication  & $8.59$s           & $27.41$s           & $150.68$s           & $490.72$s \\
					other operations       & $0.88$s           & $1.04$s            & $5.50s$             & $3.00$s  \\
				\hline
			\end{tabular}}		
		\end{center}
		\caption{Time in seconds needed by Advanpix (line 2) and Algorithm~\ref{alg:mixedprec} (lines 3--5).}
		\label{tbl:example2-timings}
	\end{table}

\end{example}

\begin{example}
	\label{example3-failure}
This example aims to provide insight into how Algorithm~\ref{alg:mixedprec}  fails to converge in exceptional situations, when eigenvalues are poorly separated or, equivalently, the triangular matrix equation~\eqref{eq:linearsystem} is very ill conditioned. To illustrate how such a situation affects Algorithm~\ref{alg:mixedprec}, we generate the matrix $A = XDX^{-1}$ of size $150$. Here, $X$ is a random matrix of condition $10^5$. The matrix $D$ is a diagonal matrix with diagonal elements chosen uniformly at random between $-10$ and $10$, with the exception of two clusters of size $10$. In each cluster, the cluster center is chosen randomly, and the cluster elements are perturbed randomly by at most $10^{-5}$ from the center. 

We execute Algorithm~\ref{alg:mixedprec} for computing the complex Schur form using double-quadruple precision.
Figure \ref{fig:example3} shows the absolute values of the computed matrix $L$ during the second and sixth iteration of the algorithm. We note that, initially, all the entries of $L$ are quite small with the exception of those that correspond to the two clusters. In later iterations, these elements have polluted the entire lower triangular part of $L$. In the eight iteration, the algorithm fails as the matrix $L$ contains NaNs.

Note that if this example is slightly softened (e.g., lowering the condition of $X$ to $10^4$ or the increasing the cluster radius to $10^{-4}$) then Algorithm~\ref{alg:mixedprec} converges in $6$ iterations and yields a Schur decomposition in quadruple precision.

	\begin{figure}
		\centering
		\includegraphics[width=0.49\textwidth]{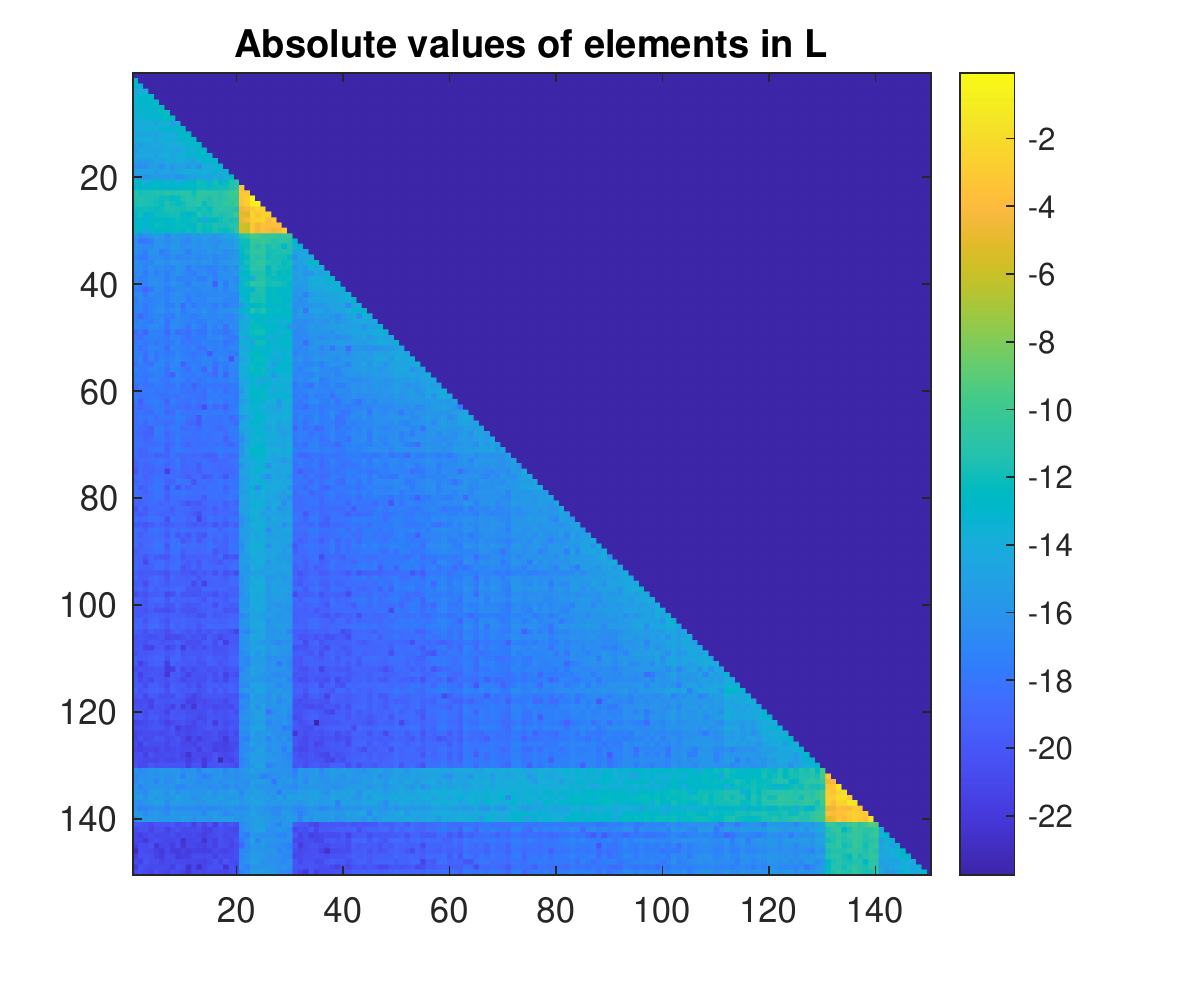}
		\hfill
		\includegraphics[width=0.49\textwidth]{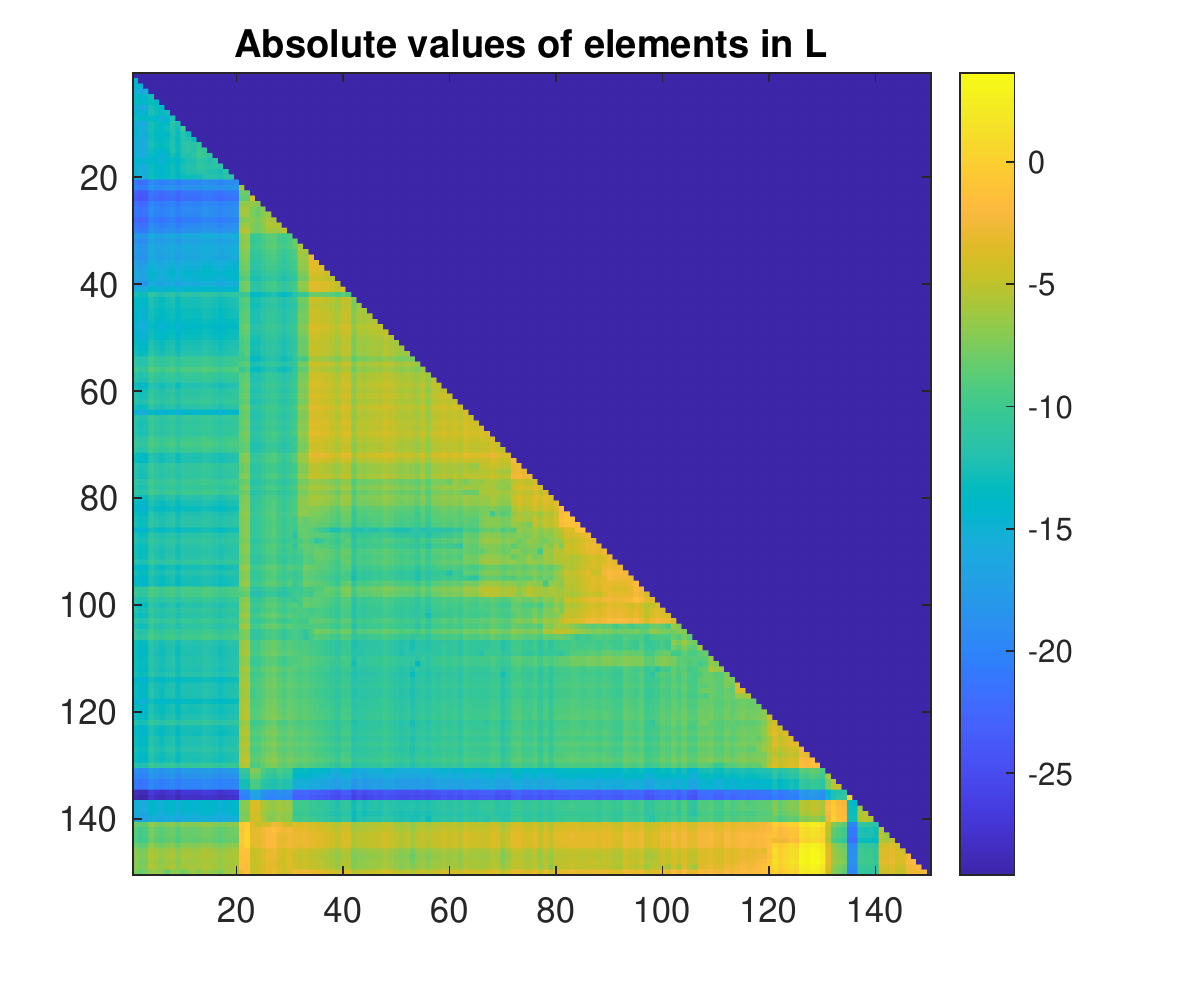}
		\caption{Absolute values of entries of the matrix $L$ obtained from solving the triangular matrix equation~\eqref{eq:linearsystem} within Algorithm~\ref{alg:mixedprec} applied to the matrix from Example \ref{example3-failure}. Left: Second iteration of Algorithm~\ref{alg:mixedprec}. Right: Sixth iteration. Color indicates base-10 logarithm of each matrix element's absolute value. }
		\label{fig:example3}
	
	\end{figure}
\end{example}

\begin{example} \label{example:benchmarks}
	Finally, we report on the performance of Algorithm~\ref{alg:mixedprec} for a number of well-known eigenvalue benchmark examples from MatrixMarket 
	\cite{BaiDDD97}. Table \ref{tbl:example4-matrixmarket} details the timing, the errors in the computed factors, and the number of iterations needed for each of these test cases. Note that the matrix size is indicated in the matrix name.

	Algorithm~\ref{alg:mixedprec} fails to converge for matrices denoted by asterisks, due to the reasons explained in Example~\ref{example3-failure}. However, there is a simple trick to avoid the appearance of large numbers in the computed factor $L$: during its computation, each value larger than, e.g., $10^{-5}$ in absolute value is immediately set to zero. Currently, there is little theoretical support that such a modification will lead to convergence of Algorithm~\ref{alg:mixedprec}---in fact, it still fails for the matrix from Example \ref{example3-failure}. However, for all matrices reported in Table \ref{tbl:example4-matrixmarket} this trick successfully resolves convergence problems. We implemented this trick only for the complex Schur form and the obtained results are shown in Table~\ref{tbl:example4-matrixmarket-withtrick}.
	
	\begin{table}
		\begin{center}
			{\tiny
			\begin{tabular}{|c||c|c|c|c|c|c|}
				\hline
					matrix & \makecell{time \\ (quad Schur)} & \#matmuls & \makecell{time outside \\ of matmuls} & $\|I-Q^\ast Q\|_F $ & $\frac{\|\text{stril}(Q^\ast A Q)\|_F}{\|A\|_F}$ & \#iterations \\
				\hline
				olm100    & 0.17s (0.12s)     & 12 & 0.04s & 2.73e-32 & 1.16e-33 & 3 \\
				tols90    & 0.09s (0.09s)     & 12 & 0.02s & 2.75e-32 & 2.72e-34 & 3 \\
				tub100    & 0.11s (0.15s)     & 12 & 0.03s & 2.69e-32 & 2.92e-33 & 3 \\
				bwm200    & 0.21s (0.96s)     & 12 & 0.06s & 4.26e-32 & 3.22e-33 & 3 \\
				olm500    & 1.25s (7.00s)     & 12 & 0.17s & 6.48e-32 & 1.56e-33 & 3 \\
				olm1000   & 13.61s (50.56s)   & 16 & 0.83s & 8.07e-32 & 8.02e-34 & 4 \\
				tub1000   & 8.48s (90.58s)    & 12 & 0.87s & 7.51e-32 & 2.64e-33 & 3 \\
				bwm2000   & 60.15s (737.59s)  & 12 & 4.24s & 1.28e-31 & 3.77e-33 & 3 \\
				dw2048    & 98.50s (982.43s)  & 12 & 4.07s & 1.27e-31 & 3.48e-33 & 3 \\ \hline
				rdb200*   & \multicolumn{6}{c|}{did not converge} \\
				rdb200l*  & \multicolumn{6}{c|}{did not converge} \\
				tols1090* & \multicolumn{6}{c|}{did not converge} \\
				rdb1250*  & \multicolumn{6}{c|}{did not converge} \\ \hline
			\end{tabular}}		
		\end{center}
		\caption{Real Schur decomposition of various matrices from MatrixMarket; double precision is refined to quadruple precision. The second column shows the total time in seconds for Algorithm~\ref{alg:mixedprec}; the time needed by Advanpix's {\tt schur} is shown in parentheses. The third column shows the total number of quadruple-precision matrix-matrix multiplications needed by Algorithm~\ref{alg:mixedprec}, and the fourth column shows the time the algorithm spent outside these operations. The next two columns show the errors in the computed refined factors $Q$ and $T$. The last column shows the number of refinement iterations. Algorithm~\ref{alg:mixedprec} fails to converge for the matrices denoted by asterisk.}
		\label{tbl:example4-matrixmarket}
	\end{table}

	\begin{table}
		\begin{center}
			{\tiny
			\begin{tabular}{|c||c|c|c|c|c|c|}
				\hline
					matrix & \makecell{time \\ (quad Schur)} & \#matmuls & \makecell{time outside \\ of matmuls} & $\|I-Q^\ast Q\|_F $ & $\frac{\|\text{stril}(Q^\ast A Q)\|_F}{\|A\|_F}$ & \#iterations \\
				\hline
				rdb200   & 0.91s (1.86s)     & 16 & 0.07  & 3.99e-32 & 3.36e-33 & 4 \\
				rdb200l  & 0.88s (2.27s)     & 16 & 0.07  & 3.99e-32 & 2.71e-33 & 4 \\
				tols1090 & 135.45s (67.77s)  & 20 & 1.95  & 8.93e-32 & 2.22e-35 & 5 \\
				rdb1250  & 143.98s (528.12s) & 16 & 1.90  & 9.79e-32 & 3.58e-33 & 4 \\ \hline
			\end{tabular}}		
		\end{center}
		\caption{The trick described in Example~\ref{example:benchmarks} resolves the convergence failures reported in Table~\ref{tbl:example4-matrixmarket} when refining the complex Schur decomposition.}
		\label{tbl:example4-matrixmarket-withtrick}
	\end{table}	
\end{example}

\section{Conclusions}

In this work, we have developed iterative algorithms for refining approximate Schur decompositions that exhibit rapid convergence, in theory and in practice. Using the Newton-Schulz iteration for orthogonalization yields an algorithm that carries out most operations in terms of matrix-matrix multiplications, allowing for a simple and efficient implementation. In particular, when refining a double precision decomposition to high precision this allows to leverage the Ozaki scheme and attain significant speedup over existing implementations of algorithms for computing high-precision Schur decompositions.

A number of points deserve further investigation. It is not unlikely that a suitable extension of the modifications described for the symmetric case in~\cite{Ogita2019} will address 
the convergence failures observed in Examples~\ref{example3-failure} and~\ref{example:benchmarks}. Also, it would be valuable to further study possibilities for merging the improvement of orthogonality and triangularity in a single step, as in~\cite{Ogita18}, and avoiding the need for the initial orthogonalization in Algorithm~\ref{alg:mixedprec}. Finally, to attain very high precision it would certainly be beneficial to study the effective use of more than two levels of precision.

\begin{paragraph}{Acknowledgments.}
The authors thank Takeshi Ogita for providing the Matlab toolbox \code{acc} based on \cite{Ozaki12}. They are also grateful to Nicolas Boumal and Christian Lubich for discussions related to Remark~\ref{remark:landing}.
\end{paragraph}

\bibliographystyle{plain}
\bibliography{all}

\end{document}